# Security Constrained AC Transmission Network Expansion Planning

Soumya Das [1*], Ashu Verma [2], P. R. Bijwe [2]

[1] Centre for Energy Studies, Indian Institute of Technology Delhi, New Delhi, India
[2] Department of Electrical Engineering, Indian Institute of Technology Delhi, New Delhi, India
*soumya.das@ces.iitd.ac.in

**Abstract:** Modern transmission network expansion planning (TNEP) is carried out with AC network model, which is able to handle voltage and voltage stability constraints. However, such a model requires optimization with iterative AC power flow model, which is computationally so demanding that most of the researchers have ignored the vital (N-1) security constraints. Therefore, the objective of this research work is to develop an efficient, two stage optimization strategy for solving this problem. In the first stage, a DC expansion planning problem is solved which provides an initial guess as well as some very good heuristics to reduce the number of power flow solutions for the second stage of AC transmission and reactive expansion planning. A modified artificial bee colony (MABC) algorithm is used to solve the resulting optimization problem. Static AC TNEP results for Garver 6 bus, IEEE 24 bus and IEEE 118 bus test systems have been obtained with the proposed and rigorous approaches and wherever possible, compared with similar results reported in literature to demonstrate the benefits of the proposed method. Also, multi-stage dynamic AC TNEP for the Garver 6 bus system is solved to show the applicability of the methodology to such problems.

## Nomenclature

| Symbol | Description |
|---|---|
| $v$ | total investment cost |
| $v_0$ | cost of line addition |
| $v_1$ | cost of additional reactive power sources |
| $\Omega$ | set of all right of ways or power corridors, with each corridor consisting of identical lines |
| $l$ | right of way (ROW) between buses $i$ and $j$ consisting of identical lines |
| $C_t^l$ | cost of adding $t^{th}$ circuit in the $l^{th}$ right of way |
| $n_l$ | number of circuits added in $l^{th}$ right of way |
| $N_{pq}$ | load buses of the system |
| $N_{pqbus}$ | set of all load buses |
| $N_{pvbus}$ | set of all generator buses including slack |
| $c_{0N_{pq}}$ | fixed cost of adding a reactive power source at the $N_{pq}{}^{th}$ load bus |
| $c_{1N_{pq}}$ | variable cost of reactive power source at the $N_{pq}{}^{th}$ load bus |
| $q_{reacN_{pq}}$ | size of reactive power source at $N_{pq}{}^{th}$ load bus |
| $u_{N_{pq}}$ | binary decision variable for determining whether to install a reactive source at a load bus or not |
| $k$ | status of line outage; $k = 0$ denotes the base case |
| $\boldsymbol{P}^k$ | vector of real power injection at the buses |
| $\boldsymbol{P_G}^k$ | vector of real power generation at the buses |
| $\boldsymbol{P_D}$ | vector of real power demand at the buses |
| $\boldsymbol{Q}^k$ | vector of reactive power injection at the buses |
| $\boldsymbol{Q_G}^k$ | vector of reactive power generation at the $pv$ buses |
| $\boldsymbol{Q_D}$ | vector of reactive demand at the buses |
| $q_{reac}$ | vector of size of reactive power source at the load buses |
| $\boldsymbol{V}^k$ | vector of voltage magnitudes at the buses |
| $\boldsymbol{\theta}^k$ | vector of voltage angles at the buses |
| $\boldsymbol{n}^k$ | vector of the system circuits |
| $L$ | L-index value of the system for base case |
| $n_l^0$ | number of circuits present before expansion planning |
| $\bar{n}_l$ | maximum additional circuits allowed for $l^{th}$ ROW or power corridor |
| $S_l^{kfrom}$ | apparent power flow through the $l^{th}$ ROW at the sending end |
| $F_t^l$ | maximum apparent power flow limit on $t^{th}$ circuit for the $l^{th}$ ROW |
| $S_l^{kto}$ | apparent power flow through the $l^{th}$ ROW at the receiving end |
| $nl$ | total number of ROWs in the system. |
| $NC$ | total number of system contingencies |
| $N_{bus}$ | set of all buses |
| $v_{dym}$ | investment cost referred to the first year |
| $t_y$ | year of planning |
| $Tyears$ | total number of years in planning horizon |
| $D_{t_y}$ | discount factor for the investment cost in a year |
| $v_{t_y}$ | investment cost in a year |
| $M$ | modified objective function |
| $\chi$ | vector of state variables |
| $\mu$ | vector of control variables |
| $\rho$ | vector of fixed variables |
| $\eta$ | weightage for the penalty functions for equality constraints |
| $E_g$ | the penalty functions for equality constraints |
| $H$ | the penalty functions for inequality constraints |
| $o$ | number of operating constraints in the problem |
| $\kappa_r$ | multiplying factor for the $r$-th operating constraint |
| $X_r^k$ | value of $r$-th control or state variable |
| $cs_N$ | colony size of the modified artificial bee colony (MABC) algorithm |
| $E_h$ | number of neighbours in the MABC algorithm |
| $iter$ | maximum iterations per trial |
| $lim$ | trial limit for generation of scout bees |
| $Dim$ | dimension size of the TNEP problems |
| $w_g$ | factor to control the effect of global optimum on a bee movement |
| $tp$ | time required per trial |
| $ff_n$ | number of fitness function evaluations required to obtain a solution |

## 1. Introduction

Transmission network expansion planning (TNEP) problem is a very well researched problem in power system planning studies. It is computationally very challenging due to its large dimensionality and mixed-integer nature. Numerous methods and techniques have been utilized for solving this problem [1]-[20] and [22]-[27]. However, most of the planning studies carried out in recent past focus on simplistic DC solution of the problem [1]-[10]. Solving TNEP with such simple model has the advantages like: 1) consideration of only real power for power equality and network constraints; 2) system representation by a reactive network, leading to a single step solution of the power flow equations; 3) no problem of non-convergence.

DC modelling provides a fast, non-iterative solution to the planning problem; however, in real world AC systems, the planning obtained by DC solution requires suitable reinforcements to maintain bus voltages and system voltage stability within limits. DC modelling has the following disadvantages:

1) In a real system, it can cause line overloading as line reactive power flows are not considered;
2) Difficulty in considering power losses at initial planning stages;
3) Flat system voltage profile of 1 p.u. is assumed, but in real world situation, voltages of load buses may fall outside the permissible limit.

All of these shortcomings of DC based planning can be addressed by AC modelling, but, this demands solution of iterative AC power flow (ACPF) at each step, thereby increasing the computational burden tremendously even for a small system. Other challenges of TNEP with AC modelling lie in considering the reactive power losses in the system and requires additional reactive power planning (RPP). If proper RPP is not done, voltage profile of the entire system may deteriorate, and in extreme cases, unavailability of sufficient reactive power can even cause voltage collapse or non-convergence of the ACPF, leading to an overall infeasible planning.

However, in current deregulated power system scenario with strict guidelines for network operations, maintaining a good voltage profile of the system and avoiding voltage collapse are major requirements of planning studies. Therefore, solving ACTNEP problem is gaining interest.

Existing literature have attempted to solve the ACTNEP problems. TNEP with full AC modelling is solved by CHA in [11]. Although CHA could not necessarily find the global optimum, it has shown the applicability of TNEP with full AC model. TNEP and RPP are solved in three stages in [12] with binary and real GA (RGA) algorithms. An initial DC planning stage is reinforced in the second stage by ACTNEP considering local supply of reactive demand. Third stage solves RPP for the new reinforced lines. Similar problem for a restructured power system is solved with PSO in [13]. Such sequential solution of TNEP and RPP is quite likely to produce inferior planning costs compared to that obtained by solving combined TNEP and RPP. Mixed integer conic programming (MICP) is used in [14], [15] for solving a linearized ACTNEP. In [14], the author has also solved multi-stage TNEP for 6 and 46 bus systems. Although the solution of 6 bus multi-stage planning is obtained very quickly, solution of 46 bus system involves huge time. Therefore, the author recommends to explore new techniques aimed at reducing the search space. In [15], authors have used rectangular co-ordinate formulation of ACPF to comply with the solution procedure. Benders decomposition is used in [16] to solve ACTNEP in three stages. However, even negative values for reactive compensation devices are reported, which may not reflect an optimum result. In [17], RGA combined with interior point method (IPM) is used for simultaneous solving of ACTNEP and RPP. Linear approximated TNEP is solved in [18] considering network contingencies. Of the total number of network contingencies, 11 most severe contingencies are considered to get the final planning results. Although the presented approach does not reach the optimal solution, the results are obtained very efficiently. An improved version of PSO is used in [19] for solving TNEP by consideration of shunt compensation in AC formulation. Here, the objective is to minimize the total cost of transmission line additions and the cost of active and reactive load shedding. In [20], a differential evolution algorithm is used for solving ACTNEP with total reactive power requirement of the system considered to be supplied by the active power generators. An effective linear approximation for full ACPF is presented in [21] for economic power dispatch problems. The proposed method, due to its simplicity, has substantial potential to be used in future power system planning studies. Also, approximated mixed-integer linear programming (MILP) models are used for ACTNEP in [22], [23] and [24]. Authors in these references propose a midway approach between a full AC model and an oversimplified DC model, with results nearer to that obtained with full ACPF modelling. One of the drawbacks reported in [23] is the increase in the number of variables with increase in system size. In [24], simultaneous solution of ACTNEP with RPP is presented under various levels of load and wind generations. The proposed method shows to be achieving superior results compared to the methods which use sequential ACTNEP and RPP. Authors in [25] have solved a mixed-integer linear approximated ACTNEP for static and multi-stage problems. Some very selective contingencies at some power corridors are simulated to obtain the final results. Similarly, in [26], a multi-objective ACTNEP is solved with evolutionary algorithm for some selective network contingencies based on the corrective control risk index value. As only some selective network contingencies are tested for the TNEP against all possible network contingencies, the resultant planning is much susceptible to produce infeasible network operations in at least a few of such omitted contingencies. A two-stage algorithm has been proposed in [27] for solution of ACTNEP with high penetration of renewable energy resources. Here, objective function considers both investment as well as operational costs. The entire planning horizon is divided in smaller blocks and sequential planning for each block is done. Linearized AC formulation with MILP is used by authors in [28] for solving static and multi-stage dynamic TNEP. Dealing with dimensionality and large computational burden is the main drawback faced by the method. Discrete evolutionary PSO (DEPSO) in [29] and high-performance hybrid genetic algorithm (HGA) is used in [30] to solve Static and multi-stage TNEP by DC formulation. Even though DC formulation is used, the authors have not considered network contingencies in their studies. In [31] authors have considered load uncertainty for solving DC modelling based dynamic TNEP.

To the best of authors' knowledge, for solving TNEP with full AC model, network contingencies have not been considered in most of the existing literature. Even with DC modelling also, in many cases, network security issues have been neglected by researchers, especially when solving multi-stage dynamic planning problems. It may be due to the tremendous computational burden faced in solving such problems. Even if some recent works do consider network contingencies, they do so with linearized approximate models. Further, all possible contingencies are never taken into account. Also, to the best of authors' knowledge, none of these papers have considered voltage stability aspects in TNEP. Voltage stability constraints, if considered appropriately in the planning stage, may lead to a robust network design with increased load margins. However, in such cases one needs to find the trade-off between stability improvement and cost of TNEP.

In view of the above, the motivation in this work is to develop suitable strategies which make it feasible to solve static, and multi-stage dynamic, security constrained (N-1 contingency) ACTNEP problems. It is achieved by solving compensation based [32] DCTNEP problem in the first stage. This requires very small fraction of the total computational effort. However, it provides vital heuristics to dramatically reduce the number of power flow solutions required in the subsequent, second stage involving ACTNEP. Both stages involve the (N-1) security constraints as well. The voltage stability constraints are added in TNEP model using L-index formulation for base case. A modified artificial bee colony (MABC) algorithm is used to solve the TNEP problems. Although there are multiple challenges in a TNEP, like, uncertainty handling, market-based requirements, etc.; this work focuses only on the basic problem to demonstrate how AC problem formulation can be used in presence of contingency constraints. As discussed earlier, this in itself is a major computational challenge for a realistically sized system. Handling other challenges along with the above requires significant amount of further research which is beyond the scope of this paper. Also, the RPP is done in this work in a simplistic way considering all load buses as candidate locations for installing additional reactive power sources. More efficient RPP techniques can be explored in future work. The main contribution of this paper lies in the development of suitable strategies for obtaining efficient solutions to the ACTNEP problems in presence of network contingencies.

The remaining part of this paper is organized as follows. Section 2 discusses the AC modelling of TNEP. Section 3 presents the heuristics applied for reduction of computational burden. Section 4 provides the simulation results and discussion for establishing the computational efficiency of the proposed approach. Section 5 draws conclusions for this work and states its future prospects.

## 2. TNEP with AC Model
### 2.1. Static TNEP

For a given loading condition, solving static ACTNEP problem involves finding the lowest total cost of line and additional reactive power sources, such that, the network equality and inequality constraints are satisfied. Each line has a power flow limit and each bus has upper and lower voltage bounds. These inequality and equality constraints are ascertained by solving ACPF. Further, each ROW can accommodate up to a certain number of lines.

L-index [33] can be considered as a fair indicator for the voltage stability of a system. Here, system voltage stability constraint in terms of L-index voltage collapse performance index (VCPI), is considered, only for the base case. The purpose of this work is to improve system stability at least in the base case, which will automatically lead to subsequent increment of stability for N-1 contingencies once the final plan is achieved. However, the developed TNEP model is general enough and if a planner wishes, the stability constraints for contingency cases may also be included without any algorithmic change.

The mathematical model for static ACTNEP with RPP as formulated in [12] can be modified for consideration of security and stability constraints as below:

Minimize:
$$v = v_0 + v_1 \qquad (1)$$
where,
$$v_0 = \sum_{l \in \Omega} \sum_{t=1}^{n_l} C_t^l \qquad (2)$$
$$v_1 = \sum_{N_{pq} \in N_{pqbus}} (c_{0N_{pq}} + c_{1N_{pq}} q_{reacN_{pq}}) u_{N_{pq}} \qquad (3)$$
such that:
$$P(V,\theta,n)^k - P_G^k + P_D = 0$$
$$Q(V,\theta,n)^k - Q_G^k + Q_D - q_{reac} = 0$$
$$P_{Gmin} \le P_G^k \le P_{Gmax}$$
$$Q_{Gmin} \le Q_G^k \le Q_{Gmax}$$
$$q_{reacmin} \le q_{reac} \le q_{reacmax}$$
$$V_{min} \le V^k \le V_{max}$$
$$L_{min} \le L \le L_{max} \qquad (4)$$

for $l \in 1,2,\ldots,nl$ and $l \ne k$,
$$S_l^{kfrom} \le \sum_{t=1}^{n_l^0 + n_l} F_t^l \qquad (5)$$
$$S_l^{kto} \le \sum_{t=1}^{n_l^0 + n_l} F_t^l \qquad (6)$$
for $l = k, k \ne 0$

$$S_l^{kfrom} \leq \sum_{t=1}^{n_l^0 + n_l - 1} F_t^l \quad (7)$$

$$S_l^{kto} \leq \sum_{t=1}^{n_l^0 + n_l - 1} F_t^l \quad (8)$$

$$0 \leq n_l \leq \bar{n}_l \quad (9)$$

$n_l \geq 0$ and integer for $l \in 1, 2, \ldots, nl$ and $l \neq k$; $(n_l^0 + n_l - 1) \geq 0$ and integer for $l = k$, $k \neq 0$. $k = 0, 1, \ldots NC$, denotes the contingency state of the system. Subscripts $min$ and $max$ respectively denote the minimum and maximum values of the associated variables.

Here, (1) represents the total investment cost which is the sum of line addition costs, represented by (2) and reactive addition cost, represented by (3). Equation (4) represents the typical hard and soft constraints of an AC optimal power flow problem. They represent real and reactive power equalities at each bus; limits on the active and reactive power generations; and, limits on the size of additional reactive power sources, voltage magnitudes and base case L-index value of the system. Equations (5) and (6) set the line power flow limits at sending and receiving ends respectively for the base case and at power corridors without any contingency. Similarly, (7) and (8) represents the line power flow limits for the corridor having a contingency and, (9) sets the limit on maximum number of circuits in a particular ROW or power corridor.

Here, $q_{reac}$ is a reactive control vector which remains fixed in base case and also after any line contingency. This is a very simplistic way of RPP, done with the objective to reduce the overall planning cost and obtain reactive power balance of the system for successful ACTNEP. Further, it should be noted that, in the TNEP model used here, whenever there is an addition of a new line with different parameters/ limits to an existing power corridor between two buses, it is considered as a separate sub-corridor between the same buses. That is, each physical power corridor between two buses is sub-divided into several separate sub-corridors, each consisting of exactly identical lines. Although such treatment increases the total number of corridors to be handled by the TNEP problem, it is done so that the power flow through the new lines can be kept within specified safe limit. Therefore, '$l$' basically denotes these sub-corridors between buses '$i$' and '$j$'.

### 2.2. Dynamic TNEP

Multi-stage dynamic TNEP not only provides which line to install, but also the optimal time for its installation over a planning horizon so that the cumulative overall cost referred to the first year becomes minimum. Cost of additional line and reactive power sources is referred to the first year with a discounted cost represented by discount factors. Extension of the static TNEP problem to dynamic is relatively straight forward, such as [14]:

Minimize:

$$v_{dym} = \sum_{t_y=1}^{Tyears} (D_{t_y} \times v_{t_y}) \quad (10)$$

Here, (10) represents the total investment cost referred to the first year which is the discounted sum of the total investment cost for each year, $v_{t_y}$, represented by (1). $Tyears$ represent the total number of years in the planning horizon and $D_{t_y}$ is the discount rate for individual year. The constraints (4) to (9) need to be satisfied for each year, with every variable augmented with subscript $t_y$ to denote the respective year. For a static problem, $Tyears$ is 1. In addition, investments done at a previous stage are always present in the later stages, such that:

$$u_{N_{pq_{t_y}}} \geq u_{N_{pq_{t_{y-1}}}} ; q_{reacN_{pq_{t_y}}} \geq q_{reacN_{pq_{t_{y-1}}}} \quad (11)$$

$$n_{l_{t_y}} \geq n_{l_{t_{y-1}}} \quad (12)$$

Multi-stage planning problems pose a tremendous challenge to the optimization algorithms as optimum planning is obtained by considering the entire planning horizon in contrast to a single year scenario in static problems. It involves a look into the future so as to obtain minimum cumulative investment cost over the planning period.

### 2.3. Algorithm used to solve TNEP

In this research work, the ACTNEP and RPP problems are solved using MABC algorithm [34], which is a general optimization algorithm developed from the artificial bee colony (ABC) algorithm [35], [36]. MABC is developed by incorporating the concept of universal gravitation [37] and global attraction [5] in the original ABC algorithm. It can be used for solving both constrained and unconstrained optimization problems. Like any other meta-heuristic technique, it can be applied directly to unconstrained problems. However, for its application to constrained optimization problems, the constrained objective functions need to be converted into unconstrained functions by incorporating suitable penalty terms corresponding to each constraint. The imposed penalties force the algorithm to direct its search in the feasible region of the search space. When there are no violations, values of the extra added functions become zero and the modified objective function becomes same as the original one.

For static ACTNEP, modified objective function can be written as:

$$M(\chi, \mu, \rho) = v(\chi, \mu, \rho) + \eta E_g(\chi, \mu, \rho) + H(\chi, \mu, \rho) \quad (13)$$

$E_g(\chi, \mu, \rho)$ are the penalty functions for equality constraints and $\eta$ defines their weightage.

$H(\chi, \mu, \rho)$ are the penalty functions for inequality constraints, defined as:

$$H(\chi, \mu, \rho) = \sum_{r=1}^{o} h_r(\chi, \mu, \rho) \quad (14)$$

where,

$$h_r(\chi, \mu, \rho) = \begin{cases} \kappa_r(|X_{rmin}| - |X_r^k|)^2 & for\ |X_r^k| < |X_{rmin}| \\ 0 & for\ |X_{rmin}| \leq |X_r^k| \leq |X_{rmax}| \\ \kappa_r(|X_r^k| - |X_{rmax}|)^2 & for\ |X_r^k| > |X_{rmax}| \end{cases} \quad (15)$$

Here, $\kappa_r$ is the multiplying factor for the $r$-th operating constraint $X_r(\chi, \mu, \rho)$, whereas $X_{rmin}$ and $X_{rmax}$ are its corresponding minimum and maximum limits. State variables ($\chi$) are: $V_i$ ($\forall i \in N_{pqbus}$) and $\theta_i$ ($\forall i \in N_{bus}, i \neq$ slack). Control variables ($\mu$) are: $P_{Gi}$ ($\forall i \in N_{pvbus}$ $i \neq$ slack), $q_{reaci}$ ($\forall i \in N_{pqbus}$), $V_i$ ($\forall i \in N_{pvbus}$) and $n_l$ ($\forall l \in \Omega$). Fixed variables ($\rho$) are: $P_{Di}(\forall i \in N_{pqbus})$, $Q_{Di}(\forall i \in N_{pqbus})$, $\theta_i(i =$slack).

It is evident from (14) and (15) that, $H$ is zero when no operating limits are violated. Also satisfying the equality constraints make $E_g = 0$.

When considering multi-stage problems, in (13), $v(\chi,\mu,\rho)$ is replaced by $v_{dym}(\chi,\mu,\rho)$; $E_g(\chi,\mu,\rho)$ and $H(\chi,\mu,\rho)$ are replaced by $E_{gdym}(\chi,\mu,\rho)$ and $H_{dym}(\chi,\mu,\rho)$, which are represented as:

$$E_{gdym}(\chi,\mu,\rho) = \sum_{t_y=1}^{Tyears} E_{g_{t_y}}(\chi,\mu,\rho) \qquad (16)$$

$$H_{dym}(\chi,\mu,\rho) = \sum_{t_y=1}^{Tyears} H_{t_y}(\chi,\mu,\rho) \qquad (17)$$

Now, a general algorithm for solving TNEP can be described as follows:

Step 1: Start the algorithm
Step 2: Read—line, branch and load data. Form the modified objective function, $M(\chi,\mu,\rho)$.
Step 3: Take as input the control parameter values of MABC like colony size $cs_N$; maximum iteration number $iter$, $lim$, $w_g$ and number of neighbours $E_h$. Determine the dimension size, $Dim$ of the problem.
Step 4: Initialize the MABC algorithm by random generation of candidate solutions within the search space.
Step 5: Set $itercycle = 1$
Step 6: For each of the solution points thus produced by each stage of the MABC algorithm, ACPF is solved. From the ACPF solutions, the value of $M(\chi,\mu,\rho)$ is calculated for each point, and then their $fitness$ values are determined. For a minimization problem, $fitness = \frac{1}{M}$.
Step 7: Store the one with the best $fitness$ value, i.e. $best\ solution$.
Step 8: Set $itercycle = itercycle + 1$
Step 9: Repeat—Step 6 to 8 until $itercycle = iter$.
Step 10: Print $best\ solution$ and the results of the TNEP.
Step 11: End the algorithm

## 3. Heuristics applied for the reduction of computational burden

Solving an ACTNEP problem with N-1 contingency constraints involves high complexity when TNEP and RPP are solved together, while also considering the voltage stability constraints. To efficiently solve such a challenging problem, several heuristics are applied to the MABC algorithm. The proposed heuristics lead to a manageable computational burden, which can be handled to get the solution of security constrained ACTNEP problem. Here, TNEP is solved in two stages. Further, few heuristics are developed to make the process more efficient in the second stage. The heuristics applied are described as follows:

### 3.1. 1st Stage: Starting the ACTNEP from an initial DC solution:

An appropriate initial assumption is required for a metaheuristic algorithm to reach an acceptable solution within a short time. Here, the solution obtained by compensation based [32] DCTNEP is used as a starting point for the ACTNEP. DCTNEP solution is obtained with negligible computational effort compared to ACTNEP, and provides a good, sub-optimal starting point for the latter.

### 3.2. 2nd Stage Heuristics:

#### 3.2.1. Reducing the number of power flow solutions:

For solving ACTNEP, huge number of combinations are required to be evaluated. This involves enormous time to search, as for each combination, it is required to solve ACPF. Time taken for obtaining the optimum solution by a metaheuristic algorithm depends on the number of times power flow is solved. Lesser this number, lesser time it will take to reach the optimum.

It has been found by experience that almost 90% of the corridors for new lines found by the DC solution are also present in AC solution. Obviously, ACTNEP requires some more lines. Thus, to effectively reduce the number of times ACPF solution is required in MABC algorithm, the total number of power corridors for ACTNEP is restricted to be within 90% and 130% of the number of corridors that is obtained by DC solution. If the number of power corridors suggested by a combination is within this range, ACPF is solved. In other cases, a suitable penalty is added.

#### 3.2.2 Checking the worthiness of a solution before actually running a power flow and evaluating its fitness function:

In a metaheuristic algorithm, large number of combinations are produced which are tested for their feasibility. However, it is observed that, most of the combinations produced at the initial stages of the algorithm get rejected after their fitness values are evaluated using ACPF. Hence, huge computational time is required to evaluate the fitness values of these infeasible combinations, thereby making the algorithm computationally inefficient.

Therefore, this heuristic focuses on determining the worthiness of a particular combination before solving ACPF. If a combination is deemed to be worthy enough, power flow is solved, otherwise, a suitable penalty is added to the objective function. Determination of worthiness depends on the cost of the suggested new lines. Here, ACPF is solved, only if the cost of new lines for the combination is less than a specified upper limit.

Exercising a tight upper bound on the criteria to discard a combination may reject too many potential combinations, leading to a population pool with a very less variation. As the effectiveness of a metaheuristic algorithm depends on variation in its candidate solutions, such tight limit may lead the same to get trapped in a local optimum. In this work, application of MABC reveals that, setting an upper limit as twice the cost obtained by the DC result gives fair balance between computational time and an acceptable final result.

This two-stage process and second stage heuristics together make the proposed method computationally very efficient by reducing the number of ACPF solutions and still providing enough flexibility to find the global optimum solution. In Section 4, it will be evident that, such actions help in obtaining feasible results for the security constrained ACTNEP problems with a drastically reduced computational burden compared to the single stage rigorous method (which does not use any of the proposed heuristics). Detailed algorithm flowchart is also presented in Fig. 1. It should be noted here that, solving ACPF in the algorithm actually involves solving it for base as well as N-1 contingency cases, for every year in the planning horizon. Cumulative penalty function values due to various limit and constraint violations are then used for the evaluation of the modified objective function.

## 4. Results and Discussion

The performance of the proposed method is evaluated by solving static ACTNEP for Garver 6 bus [12], IEEE 24 bus [11] and modified IEEE 118 bus standard systems [18]. Further, the applicability of the methodology to multi-stage dynamic planning problems has been demonstrated by solving multi-stage dynamic ACTNEP for Garver 6 bus [14] system. Static ACTNEP for 6 and 24 bus systems are solved with and without considering the security constraints for both fixed and dispatch-able generation scenarios. For 118 bus system and 6 bus dynamic system, ACTNEP is solved only for dispatch-able generation scenario, with and without considering the security constraints. Results for some cases are also obtained with rigorous method to demonstrate the potential of the proposed algorithm. All the simulations are carried out on MATLAB version R2015b running on a desktop computer with Intel (R) Core(TM) i5-4590 CPU @ 3.30 GHz and 16.0 GB RAM. The results obtained in this section, are the best among 50 trials, which is a fair number to confirm the effectiveness of any metaheuristic algorithm. Algorithm parameters for MABC as shown in Table 1 are tuned according to the procedure discussed in Section 4.1.3.

**Table 1** MABC control parameters

| Method | | $cs_N$ | $E_h$ | $lim$ | $iter$ | $w_g$ |
|---|---|---|---|---|---|---|
| Proposed | DCTNEP | 5 | 2 | 6 | 15 | 1.5 |
| | ACTNEP | 20 | 2 | 6 | 30 | 1.5 |
| Rigorous | ACTNEP | 20 | 2 | 6 | 30 | 1.5 |

In [11]-[20] and [22]-[28], ACTNEP is done with generation rescheduling. Through generation rescheduling, much lower cost of line and reactive power addition is possible, since for any network configuration (base and contingency cases), violations in line power flow and bus voltage magnitudes may be reduced by proper rescheduling of the generator power and their terminal voltages. However, in a deregulated scenario, a transmission plan must not restrict the generation dispatch, which is primarily determined by economic considerations. Therefore, in addition to the results with generation rescheduling (dispatch-able generation scenario), for comparison purposes, ACTNEP results considering fixed power generation plans are also presented.

Due to the unavailability of similar results in existing literature, results obtained by the proposed method with a fixed generation plan are compared with that of rigorous, single stage method. In all of these studies, system bus voltage magnitudes are limited within 0.95 and 1.05 p.u. for base and contingency cases. To ensure adequate voltage stability of a system, its base case L-index value is limited within a specified low bound of 0.45.

### 4.1. Garver 6 Bus System

This is a 6 bus system, having 15 power corridors with a total real power demand of 760 MW and reactive power demand of 152 MVAR. Each power corridor can accommodate a maximum of 5 lines. Resistance of each line is considered to be $1/10^{th}$ of its reactance. The electrical and network data described in [11], [12] is used for this work. Installation costs of the reactive power resources are considered as 100 US$ (fixed cost) plus 0.3 US$/kvar (variable cost) [12].

*1) TNEP without security constraints:*

Base case ACTNEP result obtained for the 6 bus system with generation rescheduling, by the proposed method is similar to that reported in literature, with line addition cost of $110 \times 10^3$ US$ and reactive addition cost of $13.67 \times 10^3$ US$, thereby providing a total TNEP cost of $123.67 \times 10^3$ US$, which is lower than that obtained in [12] by 1.18%. This is a validation of the proposed method and shows its applicability in solving ACTNEP and RPP problems.

ACTNEP results for a fixed generation plan obtained by the proposed and rigorous methods are shown in Table 2. For considering a fixed generation plan, the generations at $3^{rd}$ and $6^{th}$ buses are set at 3.22 and 2.97 p.u. respectively in accordance with the results obtained for the complete test 4 in [11]. Using the proposed method resulted in computational reductions of 95.83% compared to the rigorous method.

**Table 2** ACTNEP results of Garver 6 bus system for base case (with fixed generation)

| Plan: Base case 1 | Proposed Method | | Rigorous Method | |
|---|---|---|---|---|
| New lines Constructed | $n_{2-6} = 1$; $n_{3-5} = 1$; $n_{4-6} = 2$ | | $n_{2-6} = 1$; $n_{3-5} = 1$; $n_{4-6} = 2$ | |
| No. of New Lines | 4 | | 4 | |
| Additional Reactive power sources (p.u.) | Bus 2: | 0.7373 | Bus 2: | 0.7972 |
| | Bus 4: | 0.0000 | Bus 4: | 0.0000 |
| | Bus 5: | 0.0000 | Bus 5: | 0.0000 |
| $v_0$ (x $10^3$ US$) | **110.0000** | | 110.0000 | |
| $v_1$ (x $10^3$ US$) | 22.2190 | | 24.0160 | |
| $v$ (x $10^3$ US$) | **132.2190** | | 134.0160 | |
| $L$ | 0.2206 | | 0.2172 | |
| $tp$ | 5.82 secs | | 96.78 secs | |
| $ff_n$ | 1307 | | 31361 | |
| % Reduction in Computational Burden by Proposed Method | **95.83** | | | |

*2) TNEP with security constraints:*

The best solutions obtained by the proposed and rigorous methods for ACTNEP with N-1 contingency, considering both dispatch-able and fixed generation scenarios, are given in Table 3. In the dispatch-able generation scenario, active power generations and *pv* bus voltage magnitudes are modified during each network contingency to alleviate line power flow violations. Whereas, for fixed generation scenario, these remain fixed in all network configurations. The fixed generation values are similar to that used in the base case solution. Unavailability of similar results in existing literature prevents any comparison with other methods. However, the effectiveness of the proposed method is demonstrated by comparing the results with the ones obtained by the rigorous method. From Table 3, it can be observed that with the proposed method, the overall computational burden is reduced by more than 94% compared to the rigorous method.

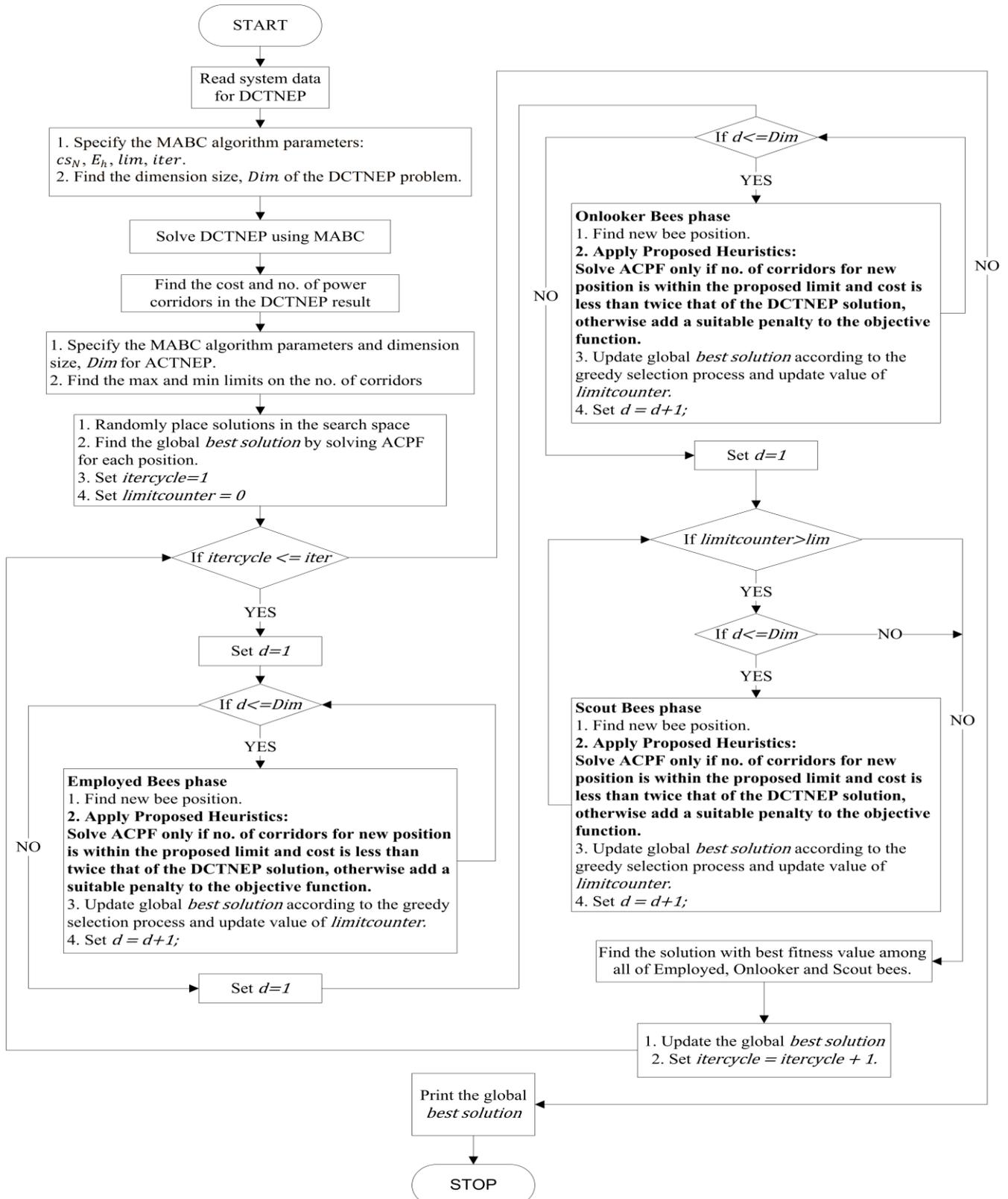

*Fig. 1.* *Flow Chart of MABC algorithm applied to the ACTNEP problem*

**Table 3** ACTNEP results obtained with the proposed and rigorous methods for Garver 6 bus system with N-1 contingency

| Plan: Base case 1 | Dispatch-able Generation | | Fixed Generation Plan | |
|---|---|---|---|---|
| | **Proposed Method** | **Rigorous Method** | **Proposed Method** | **Rigorous Method** |
| New lines Constructed | $n_{2-6}=2$; $n_{3-5}=2$; $n_{4-6}=2$ | $n_{2-6}=2$; $n_{3-5}=2$; $n_{4-6}=2$ | $n_{2-6}=2$; $n_{3-5}=2$; $n_{4-6}=2$ | $n_{2-6}=2$; $n_{3-5}=2$; $n_{4-6}=2$ |
| No. of New Lines | 6 | 6 | 6 | 6 |
| Additional Reactive power sources (p.u.) | Bus 2: 0.6369<br>Bus 4: 0.1518<br>Bus 5: 0.0000 | Bus 2: 0.6961<br>Bus 4: 0.1145<br>Bus 5: 0.0000 | Bus 2: 0.3968<br>Bus 4: 0.5934<br>Bus 5: 0.1614 | Bus 2: 0.4267<br>Bus 4: 0.5690<br>Bus 5: 0.2161 |
| $v_0$ (x $10^3$ US\$) | **160.0000** | 160.0000 | 160.0000 | 160.0000 |
| $v_1$ (x $10^3$ US\$) | 23.8589 | 24.5178 | 34.8480 | 36.6540 |
| $v$ (x $10^3$ US\$) | **183.8589** | 184.5178 | 194.8480 | 196.6540 |
| $L$ | 0.1920 | 0.1917 | 0.1874 | 0.1910 |
| $tp$ | 59.63 secs | 1389.55 secs | 40.41 secs | 720.20 secs |
| $ff_n$ | 3 046 | 72 403 | 2 358 | 43 317 |
| % Reduction in Computational Burden by Proposed Method | **95.79** | | 94.56 | |

From both the Tables 2 and 3 it can be noticed that, the proposed two stage method gives same line addition costs as the rigorous method, but, the reactive addition costs are lower than those with the latter. Such behaviour can be attributed to the fact that, in the proposed method, heuristics are applied only for line additions, whereas for reactive source additions, the search process is same as in the rigorous method. This allows the proposed algorithm to focus more on optimizing reactive source additions, thereby producing a lower reactive addition cost. However, if the values of $cs_N$ and $iter$ are increased suitably for the rigorous method, it is able to get marginally better results than the proposed method, but with huge increase of computational burdens, as shown in Table 4.

**Table 4** ACTNEP results of Garver 6 bus system with rigorous method ($cs_N$ = 50; $iter$ = 100)

| Plan: Base case 1 | Dispatch-able Generation | Fixed Generation Plan | |
|---|---|---|---|
| | Contingency Case | Base Case | Contingency Case |
| No. of New Lines | 6 | 4 | 6 |
| $v_0$ (x $10^3$ US\$) | **160.0000** | 110.0000 | 160.0000 |
| $v_1$ (x $10^3$ US\$) | 22.5573 | 21.2080 | 33.2850 |
| $v$ (x $10^3$ US\$) | **182.5573** | 131.2080 | 193.2850 |
| $L$ | 0.1858 | 0.2230 | 0.1847 |
| **COMPARED TO THE PROPOSED METHOD** | | | |
| % Decrease in overall cost | **0.71** | 0.76 | 0.80 |
| % Increase in computational burden | **12748.69** | 12363.81 | 13089.23 |

Cost convergence curves obtained by the proposed method has been provided in Figs 2a, 2b and 2c. These show that the proposed methodology obtains the optimum solution extremely fast within a very few iterations.

Table 4 also shows that, for rigorous method, with increased values of $cs_N$ and $iter$, line addition costs are similar to that obtained by the proposed method, while, total cost reductions are only due to reduction of the costs for additional reactive sources, which are fractions of the total costs. Thus, it is seen that the use of the proposed method is substantially advantageous over the rigorous method.

*3) Parameter Tuning for the MABC Algorithm:*

In any metaheuristic algorithm, parameter tuning is an important aspect that determines its efficiency to reach the global optimum solution. Parameter tuning for the MABC algorithm used in this study has been done with the target of efficiently obtaining ACTNEP solutions. For tuning of various parameters of the algorithm, solution of ACTNEP by the rigorous method, for Garver 6 bus system considering the security constraints (with dispatch-able generation) has been used.

As metaheuristic algorithm is a population dependant algorithm, a good variance in the population is needed for reaching the optimal solution. Hence, a value of the tuning parameter which exhibits a large variance of the population pool compared to the other settings is a good estimate for the optimal value of the respective parameter [34]. A few trials of the algorithm with lesser iterations per trial is enough for obtaining a good estimate as is shown in Tables 5 and 6. The control parameters of MABC are: $cs_N$, $E_h$, $lim$, $iter$ and $w_g$. However, according to [35], ABC does not show strong dependency on $cs_N$. As MABC is based on ABC, it also follows the same and the value of $cs_N$ is kept fixed at 20. Further, as will be evident from the cost convergence curves in Fig. 2 and Fig. 3, the algorithm reaches optimum solution within about 25 iterations. Therefore, the value of $iter$ is set at 30 for all the cases. The value of $w_g$ is set according to that in [5].

Therefore only $E_h$, and $lim$ values need to be properly tuned. The values of these parameters are tuned by varying one parameter at a time while keeping the other fixed. It can be observed from Tables 5 and 6 that the values which provide the highest variance also produces the least minimum, maximum and mean costs. Therefore, the values producing the highest variance is a very good estimate for optimal parameter setting and the same values have been used in all the studies conducted in this work.

**Table 5** Effect of $E_h$ on security-constrained transmission network expansion planning results for Garver 6 bus system for dispatch-able generation (with 5 trials and $lim = 6$)

| $E_h$ | | 1 | **2** | 3 | 4 | 5 | 6 |
|---|---|---|---|---|---|---|---|
| Variance of population pool | 1st trial | 14.5774 | **21.9541** | 11.4406 | 12.6303 | 12.8191 | 5.4663 |
| | 2nd trial | 15.8672 | **19.6322** | 22.0083 | 15.7411 | 11.5119 | 17.3769 |
| | 3rd trial | 17.4637 | **25.2357** | 17.8252 | 12.1962 | 6.8445 | 6.1182 |
| | 4th trial | 13.9634 | **28.4161** | 9.3749 | 10.7479 | 3.0086 | 11.8033 |
| | 5th trial | 15.2367 | **32.8932** | 14.723 | 21.3662 | 14.7965 | 5.6145 |
| Minimum cost (x $10^3$ US$) | | 266.74091 | **190.9324** | 211.9057 | 230.1339 | 215.4017 | 208.4409 |
| Maximum cost (x $10^3$ US$ | | 437.5399 | **225.3345** | 292.0389 | 239.5014 | 274.0244 | 280.1178 |
| Mean cost (x $10^3$ US$) | | 383.0316 | **208.2908** | 231.7874 | 232.9047 | 250.5009 | 243.3136 |
| Standard Deviation (x $10^3$ US$) | | 60.2297 | **13.8631** | 30.4091 | 3.4156 | 23.0190 | 26.4749 |

For the DC first stage in the proposed two-stage method, lower values of $cs_N$ and $iter$ is used only because of the fact that DC results are very easy to obtain and a final optimal DC solution is not a very essential requirement for the proposed methodology.

In relation to the setting of penalty terms for the algorithm, these are considered in accordance to the order of the objective function value. Weightages for the penalty functions are set such that a violation of a constraint produces a penalty value which is no more than 10 times the value of the original objective function.

### 4.2. IEEE 24 Bus System

This system has 24 buses, 41 power corridors and a total power demand of 8550 MW and 1740 MVAR. Electrical and circuit data described in [11] are used for this work. Installation costs of the reactive power sources are 1000 US$ as fixed cost and 3 US$/kvar as variable cost [12].

*1) TNEP without security constraints:*

ACTNEP with RPP solved for dispatch-able generation with the proposed methodology produces line addition cost of 48x$10^6$ US$, and reactive addition cost of 5.77x$10^6$ US$, producing a total cost of 53.77x$10^6$ US$ which is very close to that obtained in [12]. The same line addition cost as in [12] shows the applicability of the proposed method for the 24 bus system also.

For the solution of ACTNEP with a fixed generation scenario, generation plan G1 is used here in accordance with [1]. As similar results are unavailable in literature, effectiveness of the proposed approach is demonstrated by comparing its results with that obtained by the rigorous method. Table 7 compares the results for generation plan G1 as obtained by the proposed and rigorous methods. It shows that, the proposed method is able to get the same line addition costs as obtained by the rigorous method, while the computational burden is reduced by 85.28%. The overall investment cost obtained by the proposed method is also lower than the rigorous method due to the reason stated before. In the same way as in the 6 bus system, here also, the rigorous method obtains marginally better results if its $cs_N$ and $iter$ values are increased suitably. However, this results in tremendous increase of computational burden.

*2) TNEP with security constraints:*

In this case, due to huge computational burden involved in solving the problem by rigorous method, it has only been solved by the proposed method. Similar to 6 bus system, here also, for dispatch-able generation scenario, active power generations and $pv$ bus voltage magnitudes are modified in the contingency cases. However, in fixed generation scenario, these remain fixed for all network configurations. Table 8 provides the results of N-1 contingency constrained ACTNEP, for dispatch-able and fixed generation plan G1. It can be observed that, line costs obtained for the fixed generation scenario has increased from the dispatch-able generation scenario by 31.9%, while the overall cost has increased by 31.7%.

Table 8 also shows that, by the application of the proposed method, it is possible to get the results for N-1 contingency constrained ACTNEP problems for a moderately sized system with a reasonable computational burden. Also, the base case L-index values obtained by the proposed method are very low ($< 0.35$). If further reduction in L-index value is desired, the line addition costs become higher. For a tighter limit on the L-index value ($\leq 0.25$), the results are shown for generation plan G1 (with N-1 contingency) in Table 9. When compared with Table 8, it can be observed that restricting the L-index of the system to such a low value increases the line addition cost by 7.81%, at the cost of obtaining a greater system stability. There is also a 29.06% decrease in the cost of additional reactive sources, thereby resulting in an increase of overall planning cost by 7.73%. Therefore, by the application of the proposed method, it is possible, in a very efficient way, to choose the right network expansion plan based on the intended voltage stability. Cost convergence curves related to Table 8 are provided in Figs 2d and 3a, which depict fast convergence to final solution.

**Table 6** Effect of $lim$ on security-constrained transmission network expansion planning results for Garver 6 bus system for dispatch-able generation (with 5 trials and $E_h = 2$)

| $lim$ | | 3 | 5 | **6** | 10 | 15 | 20 |
|---|---|---|---|---|---|---|---|
| Variance of population pool | 1st trial | 12.1982 | 9.4049 | **21.9541** | 15.7524 | 4.9138 | 8.7399 |
| | 2nd trial | 16.3758 | 18.2682 | **19.6322** | 8.2873 | 18.0508 | 7.6220 |
| | 3rd trial | 3.2159 | 5.9256 | **25.2357** | 7.4625 | 20.8054 | 12.0399 |
| | 4th trial | 21.0992 | 11.6969 | **28.4161** | 6.4278 | 29.1244 | 9.0382 |
| | 5th trial | 15.9837 | 3.8596 | **32.8932** | 4.7818 | 16.3749 | 16.273 |
| Minimum cost (x $10^3$ US$) | | 208.67396 | 228.0004 | **190.9324** | 199.4019 | 217.2564 | 195.9161 |
| Maximum cost (x $10^3$ US$ | | 291.7611 | 238.9199 | **225.3345** | 228.7988 | 205.2706 | 218.8106 |
| Mean cost (x $10^3$ US$) | | 241.6390 | 231.2640 | **208.2908** | 214.8046 | 211.7575 | 2.0373 |
| Standard Deviation (x $10^3$ US$) | | 35.8367 | 3.9821 | **13.8631** | 12.1472 | 3.9746 | 8.2679 |

**Table 7** ACTNEP results of IEEE 24 bus system for base case (for fixed generation plan G1 without N-1 contingencies)

| Gen. Plan | | G1 | | | | |
|---|---|---|---|---|---|---|
| | | **Proposed Method** | | **Rigorous Method** | | **Rigorous Method** ($cs_N = 50$; $iter = 100$) |
| New Lines Constructed | | $n_{6-10} = 1$; $n_{7-8} = 2$ $n_{14-16} = 1$; $n_{16-17} = 1$; $n_{17-18} = 1$ | | $n_{6-10} = 1$; $n_{7-8} = 2$ $n_{14-16} = 1$; $n_{16-17} = 1$ $n_{17-18} = 1$ | | $n_{6-10} = 1$; $n_{7-8} = 2$ $n_{14-16} = 1$; $n_{16-17} = 1$ $n_{17-18} = 1$ |
| No. of New Lines | | 6 | | 6 | | 6 |
| Additional reactive power sources (p.u.) | Bus 3: | 2.3395 | Bus 3: | 2.9981 | Bus 3: | 2.0166 |
| | Bus 4: | 0.0000 | Bus 4: | 0.0000 | Bus 4: | 0.0000 |
| | Bus 5: | 0.0000 | Bus 5: | 0.0000 | Bus 5: | 0.0000 |
| | Bus 8: | 0.0000 | Bus 8: | 0.0000 | Bus 8: | 0.0000 |
| | Bus 9: | 3.5164 | Bus 9: | 4.0425 | Bus 9: | 3.0962 |
| | Bus 10: | 0.0000 | Bus 10: | 0.0000 | Bus 10: | 0.0000 |
| | Bus 11: | 0.0000 | Bus 11: | 0.0000 | Bus 11: | 0.0000 |
| | Bus 12: | 0.0000 | Bus 12: | 0.0000 | Bus 12: | 0.0000 |
| | Bus 17: | 0.0000 | Bus 17: | 0.0000 | Bus 17: | 0.0000 |
| | Bus 19: | 0.0000 | Bus 19: | 0.0000 | Bus 19: | 0.0000 |
| | Bus 20: | 0.0000 | Bus 20: | 0.0000 | Bus 20: | 0.0000 |
| | Bus 24: | 1.3012 | Bus 24: | 1.2812 | Bus 24: | 1.0723 |
| $v_0$ (x $10^6$ US$) | | **158.0000** | | 158.0000 | | 158.0000 |
| $v_1$ (x $10^6$ US$) | | 2.1501 | | 2.4995 | | 1.8585 |
| $v$ (x $10^6$ US$) | | **160.1501** | | 160.4995 | | 159.8585 |
| L | | 0.4075 | | 0.3905 | | 0.4179 |
| tp | | 285.71 secs | | 2046.53 secs | | 4.35 hrs |
| $ff_n$ | | 19 449 | | 132 096 | | 978 487 |
| % Reduction in Computational Burden by Proposed Method | | **85.28** | | | | |

**Table 8** ACTNEP results obtained with the proposed method for IEEE 24 bus system with N-1 contingency

| | Dispatch-able Generation Scenario | | Fixed Generation Scenario (Plan G1) | |
|---|---|---|---|---|
| New lines Constructed | $n_{1-5} = 1$; $n_{3-9} = 1$; $n_{4-9} = 1$; $n_{6-10} = 2$; $n_{7-8} = 3$; $n_{10-11} = 1$; $n_{11-13} = 1$; $n_{14-16} = 1$; $n_{14-23} = 1$; $n_{17-22} = 1$; $n_{20-23} = 1$; | | $n_{1-5} = 1$; $n_{3-24} = 1$; $n_{4-9} = 1$; $n_{6-10} = 2$; $n_{7-8} = 3$; $n_{10-12} = 1$; $n_{12-13} = 1$; $n_{14-16} = 2$; $n_{15-21} = 1$; $n_{15-24} = 1$; $n_{16-17} = 2$; $n_{16-19} = 1$; $n_{17-18} = 2$; $n_{21-22} = 1$ | |
| No. of New Lines | 14 | | 20 | |
| Additional Reactive power sources (p.u.) | Bus 3: | 3.1068 | Bus 3: | 1.4472 |
| | Bus 4: | 0.8562 | Bus 4: | 0.7581 |
| | Bus 5: | 0.6167 | Bus 5: | 0.0000 |
| | Bus 8: | 0.0881 | Bus 8: | 0.0000 |
| | Bus 9: | 2.6491 | Bus 9: | 1.5912 |
| | Bus 10: | 0.2599 | Bus 10: | 0.0000 |
| | Bus 11: | 0.0000 | Bus 11: | 0.0000 |
| | Bus 12: | 0.0000 | Bus 12: | 0.0000 |
| | Bus 17: | 0.0000 | Bus 17: | 0.0000 |
| | Bus 19: | 0.0191 | Bus 19: | 0.0000 |
| | Bus 20: | 0.0000 | Bus 20: | 0.0000 |
| | Bus 24: | 0.4602 | Bus 24: | 1.9945 |
| $v_0$ (x $10^6$ US$) | **592.0000** | | 781.0000 | |
| $v_1$ (x $10^6$ US$) | 2.4248 | | 1.7413 | |
| $v$ (x $10^6$ US$) | **594.4248** | | 782.7413 | |
| L | 0.3498 | | 0.2741 | |
| tp | 6841.62 secs | | 5934.44 secs | |
| $ff_n$ | 26 482 | | 23 341 | |

**Table 9** ACTNEP results of IEEE 24 bus system with N-1 contingency for gen. plan G1 ($L \leq 0.25$)

| New Lines Constructed | $n_{1-3} = 1$; $n_{1-5} = 1$; $n_{2-4} = 1$; $n_{3-24} = 1$; $n_{6-10} = 2$; $n_{7-8} = 3$; $n_{10-12} = 1$; $n_{12-13} = 1$; $n_{14-16} = 2$; $n_{15-21} = 1$; $n_{15-24} = 1$; $n_{16-17} = 2$; $n_{16-19} = 1$; $n_{17-18} = 2$; $n_{21-22} = 1$ | | |
|---|---|---|---|
| No. of new lines | 21 | | |
| $v_0$ (x$10^6$US$) | **842.0000** | $v_1$ (x $10^6$ US$) | 1.2353 |
| $v$ (x $10^6$ US$) | **843.2353** | L | **0.2320** |

### 4.3. Modified IEEE 118 Bus System

It is a large system with 118 buses, 179 different power corridors, 91 loads and 54 generators [38]. The total real power demand is 3733.07 MW and reactive power demand is 1442.98 MVAR with a total generator capacity of 7220 MW. As in [18], here also, line congestion is created by reducing the line capacities. However, for creating congestion, in this work, the line capacities are considered to be 60% of their original capacities. Due to the absence of actual data, the line construction costs as estimated in [18] are used for the TNEP studies. A limit of maximum of two new line constructions for each power corridor is set. Further, in this study, all line contingencies are considered for TNEP with security constraints against only 11 contingencies been considered in [18]. Therefore, direct comparison of results with that reported in earlier literature is not possible. Additional reactive sources are not added to the system load buses due to the absence of data for their installation costs. Hence, the minimization objective function for this system only includes the line construction costs. Results are only shown for dispatch-able generation scenario as specific generation plans are unavailable.

*1) TNEP without security constraints:*

In solving the base case ACTNEP problem for this system with dispatch-able generation scenario, three line additions: $n_{65-68} = 1$, $n_{80-99} = 1$ and $n_{94-100} = 1$ are obtained, resulting in an investment cost of 44.940x10$^6$ US\$. The rigorous and the proposed method both obtained similar results, however, computational reduction obtained by the proposed method over the rigorous method is almost 90%.

*2) TNEP with security constraints:*

Detailed results of ACTNEP for this case obtained by the proposed method are given in Table 10. Due to tremendous amount of computational burden involved in solving the problem by rigorous method, only the results obtained with the proposed method are presented. Active power generations and $pv$ bus voltages are modified during network contingencies similar to the previous cases. It can be observed from Table 10 that; the number of additional lines is 30 with an investment cost of 329.8688 x10$^6$ US\$. The results also show that the base case L-index value is very low. These results justify the applicability of the proposed method in solving ACTNEP for base and contingency cases even for a large system with hundreds of buses and power corridors. Corresponding cost convergence curve as depicted in Fig 3b also confirms the same.

### 4.4. 6 Bus Dynamic System

For solving the multi-stage dynamic TNEP problem for 6 bus system, greenfield expansion planning with same discount factors and load profiles as in [14] have been considered. The planning horizon considered is of three years. Study for only dispatch-able generation scenario is done, with the generation limits selected in accordance with the increased load profile in each year. The costs of reactive addition are the same that have been used previously in the static problem. Further, in the case considering network contingencies, $pv$ bus voltage magnitudes and active generations are treated as controllable variables in all the contingency scenarios. The results for TNEP without and with consideration of network contingencies are provided in Tables 11 and 12 respectively. It can be observed from the Table 11 that, in the case without network security constraints, the proposed methodology obtains planning cost of 239.6724x10$^3$ US\$, which is 2.27% lower than that obtained in [14] with L-index values within the prescribed limit. Further, the reduction of computational burden compared to the rigorous method is 96.61% with a 2.95% reduction in planning cost. For the case with N-1 security constraints, rigorous method is not applied due to extreme computational burden. Observation of Table 12 shows that here also, the results are obtained with a manageable computational burden with low L-index values in each year of planning. The cost convergence curves for this system obtained by the proposed method is presented in Figs. 3c and 3d. These show that even for an extremely complex

**Table 10** ACTNEP results obtained with the proposed method for IEEE 118 bus system with N-1 contingency

| | Dispatch-able Generation Scenario |
|---|---|
| New lines Constructed | $n_{8-5} = 1$; $n_{8-9} = 1$; $n_{8-30} = 1$; $n_{9-10} = 1$; $n_{12-117} = 1$; $n_{15-17} = 1$; $n_{17-113} = 1$; $n_{23-32} = 1$; $n_{30-17} = 1$; $n_{30-38} = 1$; $n_{38-37} = 1$; $n_{38-65} = 1$; $n_{64-65} = 1$; $n_{65-68} = 2$; $n_{68-81} = 1$; $n_{68-116} = 1$; $n_{71-73} = 1$; $n_{77-78} = 1$; $n_{80-99} = 2$; $n_{81-80} = 1$; $n_{85-86} = 1$; $n_{86-87} = 1$; $n_{94-95} = 1$; $n_{94-100} = 2$; $n_{95-96} = 1$; $n_{110-111} = 1$; $n_{110-112} = 1$. |
| No. of New Lines | 30 |
| $v_0$ (x 10$^6$ US\$) | **329.8688** |
| $v_1$ (x 10$^6$ US\$) | 0 |
| $v$ (x 10$^6$ US\$) | **329.8688** |
| $L$ | 0.0677 |
| $tp$ | 9.87 hrs |
| $ff_n$ | 31 073 |

**Table 11** Dynamic ACTNEP results of Garver 6 bus system for base case

| | | Proposed Method | | | | | Rigorous Method | | | | |
|---|---|---|---|---|---|---|---|---|---|---|---|
| | | Year 1 | | Year 2 | | Year 3 | | Year 1 | | Year 2 | | Year 3 |
| New lines Constructed | $n_{1-5} = 1$; $n_{2-3} = 1$; $n_{2-6} = 2$; $n_{3-5} = 2$; $n_{4-6} = 2$ | | | | $n_{2-6} = 1$; $n_{3-5} = 1$ | | $n_{1-5} = 1$; $n_{2-3} = 1$; $n_{2-6} = 2$; $n_{3-5} = 2$; $n_{4-6} = 2$ | | | | $n_{2-6} = 1$; $n_{3-5} = 1$ |
| No. of New Lines | 8 | | 0 | | 2 | | 8 | | 0 | | 2 |
| Additional Reactive power sources (p.u.) | Bus 2  0.0008 Bus 4  0.0365 Bus 5  0.0537 | | Bus 2  0.1728 Bus 4  0.1858 Bus 5  0.0752 | | Bus 2  0.0041 Bus 4  0.1967 Bus 5  0.0262 | | Bus 2  0.0004 Bus 4  0.0350 Bus 5  0.0084 | | Bus 2  0.3548 Bus 4  0.2088 Bus 5  0.1171 | | Bus 2  0.0758 Bus 4  0.2538 Bus 5  0.1279 |
| $v_{0_{t_y}}$ (x 10$^3$ US\$) | **200** | | **0** | | **50** | | **200** | | **0** | | **50** |
| $v_{1_{t_y}}$ (x 10$^3$ US\$) | **3.0300** | | **13.0140** | | **6.8100** | | **1.614** | | **20.421** | | **13.725** |
| $v_{dym}$ (x 10$^3$ US\$) | **239.6724** | | | | | | **246.961** | | | | |
| $L_{t_y}$ | 0.2713 | | 0.3217 | | 0.3608 | | 0.2869 | | 0.3214 | | 0.3592 |
| $tp$ | 405.77 secs | | | | | | 3.38 hrs | | | | |
| $ff_n$ | 3464 | | | | | | 102263 | | | | |
| % Reduction in Computational Burden by Proposed Method | | | **96.61** | | | | | | | | |

**Table 12** Dynamic ACTNEP results obtained with the proposed method for Garver 6 bus system with N-1 contingency

|  | Year 1 | | Year 2 | | Year 3 | |
|---|---|---|---|---|---|---|
| New lines Constructed | $n_{1-2}=1$; $n_{2-6}=3$; $n_{3-4}=1$; $n_{3-5}=4$; $n_{4-6}=2$ | | $n_{2-3}=3$; $n_{4-6}=1$; | | $n_{1-2}=1$ | |
| No. of New Lines | 11 | | 4 | | 1 | |
| Additional Reactive power sources (p.u.) | Bus 2 | 0.8231 | Bus 2 | 0.0000 | Bus 2 | 0.0000 |
|  | Bus 4 | 0.5620 | Bus 4 | 0.0000 | Bus 4 | 0.1862 |
|  | Bus 5 | 0.1184 | Bus 5 | 0.2769 | Bus 5 | 0.0000 |
| $v_{0_{t_y}}$ (x $10^3$ US$) | 329 | | 90 | | 40 | |
| $v_{1_{t_y}}$ (x $10^3$ US$) | 45.378 | | 8.307 | | 5.607 | |
| $v_{dym}$ (x $10^3$ US$) | 467.8439 | | | | | |
| $L_{t_y}$ | 0.1887 | | 0.1583 | | 0.1774 | |
| $tp$ | 1842.281 secs | | | | | |
| $ff_n$ | 3681 | | | | | |

multi-stage problem, the proposed methodology obtains the result within a few iterations.

## 5. Conclusion

This paper presents solution of contingency and voltage stability constrained AC transmission expansion planning problem. These vital aspects had been ignored earlier because of excessive computational burden involved. This paper has overcome the problem by developing a two stage optimization strategy. The first stage uses DCTNEP which provides an initial guess and good heuristics for the second stage which uses AC model. The heuristics drastically reduce the number of ACPFs to be performed, improving the computational efficiency of the process, thereby providing the non-linear, contingency constrained ACTNEP solutions, which were never attempted in the past.

From the case studies for the 6 and 24 bus systems it can be observed that for the base case scenarios with dispatch-able generation, the overall costs of expansion plans obtained by the proposed method are either better or very close to the results reported in literature. Also, for 118 bus system, the proposed method is able to get the optimal results with a drastically reduced computational burden. Also, when applied to multi-stage problem for the 6 bus system, it produced base case results better than that reported in literature. This essentially proves the robustness and effectiveness of the proposed method in solving both static and multi-stage dynamic TNEP problems. For base cases with fixed generation plans, the results obtained by the proposed method are shown to be better than those obtained by the rigorous single stage method with similar specifications. Further, the results are obtained with more than 85% reduction in computational burden compared to the rigorous method in all these cases. Similar reduction in computational burden is also shown to be achieved for N-1 contingency cases for the 6 bus system. The methodology proposed in this paper is also able to limit the base case L-index values for all the systems within a specified low value,

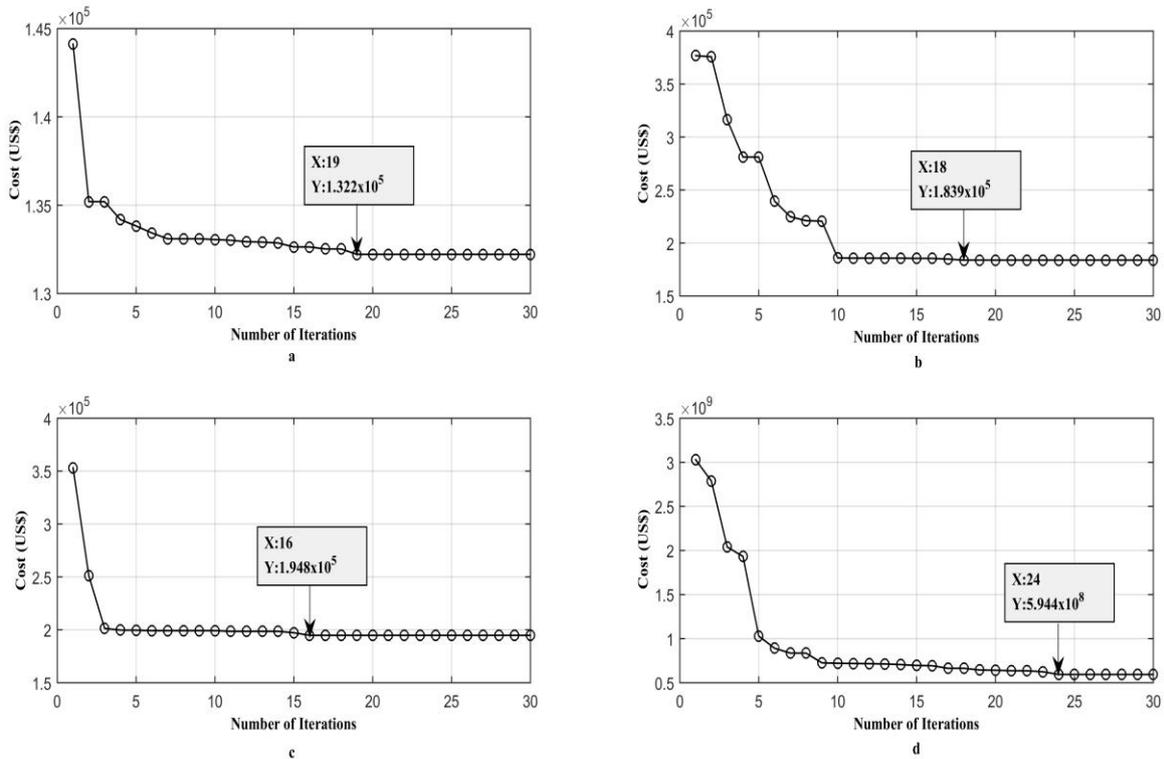

*Fig. 2.* Cost convergence curves for the static 6 bus study cases and static IEEE 24 bus system with N-1 contingency (for dispatch-able generation)
*(a)* Static 6 bus system for base case (with fixed generation), **(b)** Static 6 bus system with N-1 Contingency (for dispatch-able generation), **(c)** Static 6 bus system with N-1 Contingency (for fixed generation), **(d)** Static 24 bus system with N-1 Contingency (for dispatch-able generation)

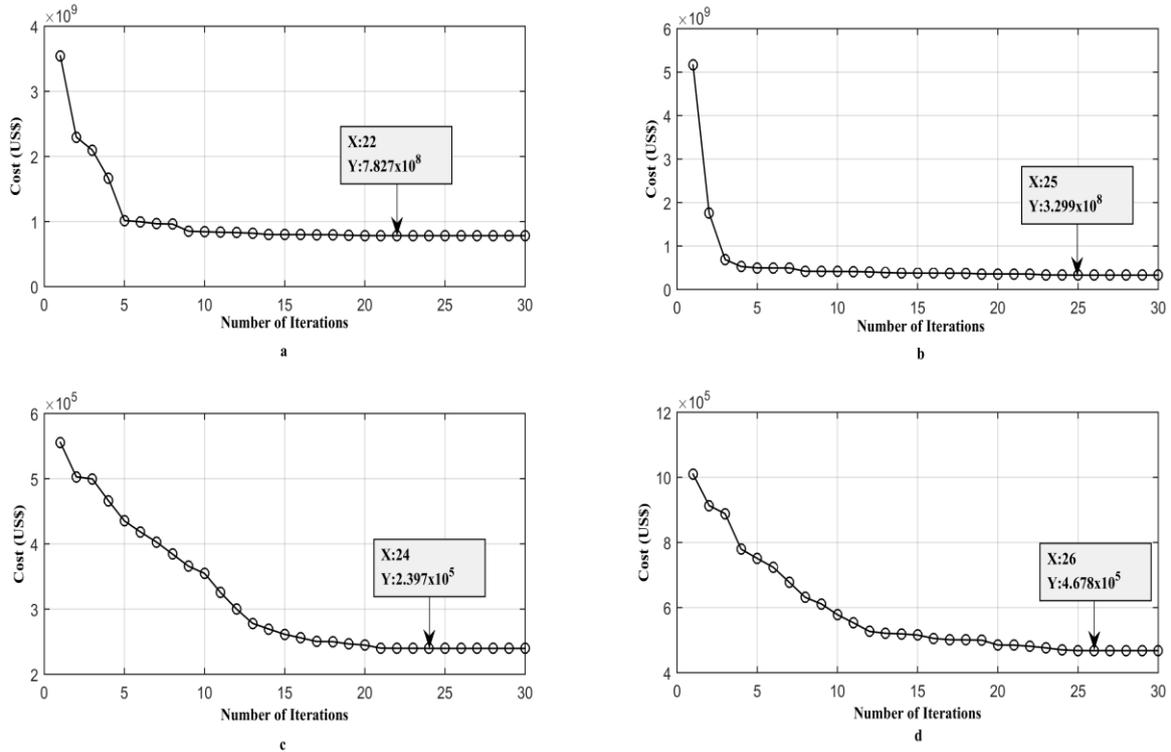

*Fig. 3.* *Cost convergence curves for the static IEEE 24 bus system with N-1 contingency (for fixed generation), static 118 bus system with N-1 contingency and dynamic 6 bus system without and with N-1 contingency*
*(a)* Static 24 bus system with N-1 Contingency (for fixed generation), **(b)** Static 118 bus system with N-1 Contingency (for dispatch-able generation), **(c)** Dynamic 6 bus system for base case (with dispatch-able generation), **(d)** Dynamic 6 bus system with N-1 Contingency (for dispatch-able generation) giving a choice to the planners to specify a suitable voltage stability criterion for their respective systems.

The cost convergence curves show that for all the cases, the proposed two-stage optimization algorithm obtains the solution within only 30 iterations. Such fast convergence to the optimum solutions for all the systems studied confirm the efficiency and efficacy of the proposed algorithm. Hence, it can be inferred that the proposed method is very adept at getting an acceptable network planning with a reasonable computational burden for base case and contingency cases.

Such reduction in computational burden by this method opens up the possibility of solving future AC transmission and reactive expansion planning problems for larger systems with even greater amount of complexity, involving load and generation uncertainties and inclusion of more complex objectives like increasing dynamic stability, optimizing market objectives, etc. which were previously not attempted due to huge computational burden.